\documentclass[%
preprint,
amsmath,amssymb,
aps,
pre,
floatfix,
]{revtex4-2}

\usepackage{graphicx}
\usepackage{multirow}
\usepackage{subcaption}
\usepackage{bm}

\newcommand{\ba}[1]{\bm{a}^{(#1)}}

\newcommand{\bd}[1]{\bm{d}^{(#1)}}
\newcommand{\bc}{\bm{c}}
\newcommand{\bdelta}{\bm{\delta}}
\newcommand{\bu}{\bm{u}}
\newcommand{\bk}{\bm{k}}
\newcommand{\bv}{\bm{v}}
\newcommand{\bx}{\bm{x}}

\newcommand{\bxi}{\bm{\xi}}
\newcommand{\bH}[1]{{\cal H}^{(#1)}}
\newcommand{\f}[1]{f^{(#1)}}
\newcommand{\pp}[2]{\frac{\partial #1}{\partial #2}}
\newcommand{\dd}[2]{\frac{d #1}{d #2}}

\newcommand{\pe}{\textsl{Pe}}
\newcommand{\re}{\textsl{Re}}
\newcommand{\pr}{\textsl{Pr}}

\begin{document}
	
\title{A multiple-relaxation-time collision model by Hermite expansion}
	
\author{Xiaowen Shan}
\email{shanxw@sustech.edu.cn}
\author{Yangyang Shi}
\author{Xuhui Li}
\affiliation{
		Guangdong Provincial Key Laboratory of Turbulence Research and Applications\\
		Shenzhen Key Laboratory of Complex Aerospace Flows\\
		Department of Mechanics and Aerospace Engineering,
		Southern University of Science and Technology, Shenzhen, 518055, China}

\begin{abstract}
The Bhatnagar-Gross-Krook (BGK) single-relaxation-time collision model for the
Boltzmann equation serves as the foundation of the lattice BGK (LBGK) method
developed in recent years.  The description of the collision as a uniform
relaxation process of the distribution function towards its equilibrium is, in
many scenarios, simplistic.  Based on a previous series of papers, we present a
collision model formulated as independent relaxations of the irreducible
components of the Hermit coefficients in the reference frame moving with the
fluid. These components, corresponding to the irreducible representation of the
rotation group, are the minimum tensor components that can be separately relaxed
without violating rotation symmetry. For the 2nd, 3rd and 4th moments
respectively, two, two and three independent relaxation rates can exist, giving
rise to the shear and bulk viscosity, thermal diffusivity and some high-order
relaxation process not explicitly manifested in the Navier-Stokes-Fourier
equations. Using the binomial transform, the Hermite coefficients are evaluated
in the absolute frame to avoid the numerical dissipation introduced by
interpolation. Extensive numerical verification is also provided.
\end{abstract}

\maketitle

\section{Introduction}

A well-known artifact of the Bhatnagar-Gross-Krook (BGK) collision
model~\cite{Bhatnagar1954} is the uniform relaxation of the distribution
function towards its equilibrium which bound all transport coefficients,
\textit{e.g.}, the shear and bulk viscosity and thermal diffusivity, to a single
relaxation time. The most noticeable manifestation is that the Prandtl number is
un-physically fixed at unity.  For kinetic theory in continuum, a couple of
remedies~\cite{Holway1966,Shakhov1972} were suggested by explicitly dialing the
Prandtl number in the equilibrium distribution. In the context of the lattice
BGK model~\cite{Chen1998a}, the multiple-relaxation-time (MRT) LB
models~\cite{d'Humieres1992} was proposed to independently relax the
eigen-states of the \textit{discrete} distribution corresponding to the
hydrodynamic moments. The details of this decomposition of the discrete
distribution is lattice-dependent.  As the velocity sets used in this class of
MRT models are insufficient for representing the third-order moments, the unity
Prandtl number was not fixed although the original MRT model was widely used for
its improved numerical stability.

In a previous series of papers~\cite{Shan2007,Li2019,Shan2019,Li2020a}, an
Hermite expansion based lattice-independent MRT collision model was developed.
The basic idea is to expand the collision term in terms of Hermite
polynomials~\cite{Grad1949a,Shan2006b} and assign a separate relaxation time to
each of the terms.  As the viscous and thermal transports are given by the
second and third moments separately, decoupling of the two were naturally
achieved. Furthermore, using sufficiently accurate quadrature rules in the
velocity space, the energy equation can be recovered.  Although yielding an
adjustable Prandtl number, the first version~\cite{Shan2007} relaxes the raw
moments and results in a Mach-number-dependent thermal diffusivity similar to
the well known ``cubic error''~\cite{Qian1993} when the thermal diffusivity is
different from the viscosity.  The second version~\cite{Shan2019,Li2019} relaxes
the central moments and restored the Galilean invariance.  Furthermore, by using
the binomial transform, the relaxation of the central moments is carried out in
the absolute reference frame without incurring the additional numerical
dissipation due to interpolation. In a more recent work~\cite{Li2020a}, the
relaxation was applied at the even finer scope of the irreducible representation
of the rotation group which is the minimal tensor space closed under spatial
rotation. For the second moment, this yields independent shear and bulk
viscosity.

The present work gives in detail a coherent presentation of the multi-relaxation
model.  In particular, the tensor decomposition beyond the second rank is given.
For the third-rank tensor one of the two possible relaxation times is identified
with the relaxation rate of the translational kinetic energy while the other,
and all three of the fourth-rank tensor, are found to have no significant effect
to the four linear hydrodynamic modes.  Similar to the BGK operator, the
first-order spac-time discretization is found to give second-order accuracy as
long as the relaxation times are shifted by one half. The paper is organized as
the following. In Sec.~\ref{sec:theory} we give the theoretical formulation of
the model. Numerical verification of the transport coefficients using the linear
hydrodynamic modes are presented in Sec.~\ref{sec:numerical}.  In
Sec.~\ref{sec:conclusion} further discussion are given.  More mathematical
details pertaining to the binomial transform between the moving and laboratory
frames are given in the Appendix.

\section{Model construction}

\label{sec:theory}

\subsection{Velocity-space discretization by Hermite expansion}

Starting from the Boltzmann-BGK equation:
\begin{equation}
	\label{eq:bgk}
	\pp{f}t + \bxi\cdot\nabla f = \Omega(f)
	\equiv -\frac{1}{\tau}\left[f-\f{0}\right],
\end{equation}
where $f(\bx, \bxi, t)$ is the single-particle distribution function, $\bxi$ the
microscopic velocity, $\bx$ and $t$ the space and time, and $\f{0}$ the
Maxwell-Boltzmann equilibrium distribution:
\begin{equation}
	\f{0} = \frac{\rho}{(2\pi\theta)^{D/2}}
	\exp\left[-\frac{(\bxi-\bu)^2}{2\theta}\right],
\end{equation}
where $\rho$, $\bu$ and $\theta$ are the fluid density, velocity and temperature
respectively, all non-dimensionalized using the scheme given in
Ref.~\cite{Shan2006b}.

The lattice-Boltzmann equation~\cite{Chen1998a} can be formulated as a velocity
discretization of Eq.~(\ref{eq:bgk})~\cite{Shan2006b}.  Expanding $f$ in Hermite
series and truncating at order $N$:
\begin{equation}
	\label{eq:hermite-expansion}
	\f{N}(\bx, \bxi, t) = \omega(\bxi)\sum_{n=0}^N\frac 1{n!}
	\ba{n}(\bx, t):\bH{n}(\bxi),
\end{equation}
where $\omega(\bxi)\equiv (2\pi)^{-D/2}\exp(-\xi^2/2)$ is the weight function,
$\bH{n}$ the $n$-th Hermite polynomial, and $\ba{n}$ the expansion coefficients
given by:
\begin{equation}
	\label{eq:an}
	\ba{n} = \int f^N(\bxi)\bH{n}(\bxi)d\bxi.
\end{equation}
Let $\{(\bxi_i, w_i), i = 1, \cdots, d\}$ be the set of abscissas and weights of
a $Q$-th degree quadrature such that the identity:
\begin{equation}
	\int\omega(\bxi)p(\bxi)d\bxi = \sum_{i=1}^dw_ip(\bxi_i),
\end{equation}
holds for all $Q$-th degree polynomial, $p(\bxi)$.  Since $f^N(\bxi) /
\omega(\bxi)$ is a polynomial of a degree $\leq N$, defining
\begin{equation}
	\label{eq:fi}
	f_i\equiv \frac{w_i\f{N}(\bxi_i)}{\omega(\bxi_i)},
\end{equation}
as long as $Q\geq 2N$, the following \textit{isomorphism} between the discrete
distribution, $f_i$, and the moment, $\ba{n}$ can be established using
Eqs.~(\ref{eq:an}) and (\ref{eq:hermite-expansion}):
\begin{equation}
	\label{eq:iso}
	\ba{n} = \sum_{i=1}^d f_i\bH{n}(\bxi_i),
	\quad\mbox{and}\quad
	f_i = w_i\sum_{n=0}^N\frac 1{n!}\ba{n}:\bH{n}(\bxi_i).
\end{equation}
On evaluating Eq.~(\ref{eq:bgk}) at $\bxi_i$, $f_i$ obeys the following lattice
BGK equation:
\begin{equation}
	\label{eq:lbgk}
	\pp{f_i}t + \bxi_i\cdot\nabla f_i = \Omega_i
	\equiv -\frac{1}{\tau}\left[f_i-\f{0}_i\right].
\end{equation}

\subsection{Chapman-Enskog asymptotic expansion}

To determine the minimum conditions for the collision term to yield correct
hydrodynamic equations, we now briefly recap how the hydrodynamic equations are
derived from Eqs.~(\ref{eq:bgk})~\cite{Huang1987}. As required by fundamental
physics, the collision term conserves mass, momentum and kinetic energy. By
taking the corresponding moments of Eqs.~(\ref{eq:bgk}), we have the
conservation equations:
\begin{subequations}
	\label{eq:cons1}
	\begin{eqnarray}
		\dd{\rho}t + \rho\nabla\cdot\bu &=& 0,\\
		\rho\dd{\bu}t + \nabla\cdot\bm{P} &=& 0,\\
		\rho\dd\epsilon t + \nabla\bu:\bm{P} + \nabla\cdot\bm{q} &=& 0,
	\end{eqnarray}
\end{subequations}
where $d/dt\equiv\partial/\partial t +\bu\cdot\nabla$ is the \textit{material
	derivative}, and
\begin{equation}
	\label{eq:pq}
	\bm{P}\equiv\int f\bc\bc d\bc,\quad\mbox{and}\quad
	\bm{q}\equiv\frac 12\int fc^2\bc d\bc,
\end{equation}
are the \textit{pressure tensor} and \textit{energy flux} respectively.
Eqs.~(\ref{eq:cons1}) must be closed by expressing $\bm{P}$ and $\bm{q}$ in
terms of $\rho$, $\bu$ and $\theta$.  At the crudest level $f$ in
Eq.~(\ref{eq:pq}) is approximated by $\f{0}$ to yield $\bm{P}^{(0)} =
\rho\theta\bdelta$ and $\bm{q}^{(0)} = 0$, which lead to the Euler's equations.
More accurate hydrodynamic equations are obtained \textit{via} the
Chapman-Enskog asymptotic calculation~\cite{Chapman1970,Huang1987}.  Let $f
\cong \f{0} + \f{1}$ where $\f{1} \ll \f{0}$ is the first approximation of the
{\em non-equilibrium} part of the distribution.  On substituting into
Eq.~(\ref{eq:bgk}) and keeping only the leading terms on both sides, we have:
\begin{equation}
	\label{eq:ck}
	\pp{\f{0}}t + \bxi\cdot\nabla\f{0} = \Omega = -\frac{\f{1}}\tau.
\end{equation}
The left-hand-side is then written in terms of $\rho$, $\bu$, $\theta$ and their
spatial derivatives with the help of Euler's equations. By taking the moments in
Eqs.~(\ref{eq:pq}) of both sides, $\bm{P}^{(1)}$ and $\bm{q}^{(1)}$ are
expressed in terms of $\rho$, $\bu$ and $\theta$. When plugged into
Eq.~(\ref{eq:cons1}), we arrive at the Navier-Stokes equations.

As shown previously~\cite{Shan2006b}, this procedure survives the velocity-space
discretization of Eqs.~(\ref{eq:iso}). Furthermore, it was
realized~\cite{Shan2019} that the same $\bm{P}^{(1)}$ and $\bm{q}^{(1)}$ with
independent proportional coefficients are obtained as long as the following
conditions are met:
\begin{equation}
	\label{eq:def}
	\int\Omega\bc^{n} d\bc = -\frac 1{\tau_n}\int\f{1}\bc^{n}d\bc,
	\quad\mbox{for}\quad n = 2, 3,
\end{equation}
where $\tau_n$ are separate relaxation times. Here both sides vanish for $n=0$
and $1$ due to conservation of mass and momentum. The BGK operator is a special
case with $\tau_2 = \tau_3 = \tau$.  A natural generalization is to require
Eq.~(\ref{eq:def}) to hold for all $n$ up to the order supported by the
underlying quadrature so that each of the central moments is independently
relaxed with relaxation time $\tau_n$.

We now formulate the collision operator in the {\em spectral space} of Hermite
polynomials. Defining $\bv\equiv (\bxi-\bu)/\sqrt{\theta}$ and the corresponding
Hermite expansion:
\begin{equation}
	f(\bv) = \omega(\bv)\sum_{n=0}^\infty\frac 1{n!}\bd{n}(\bx, t):\bH{n}(\bv),
\end{equation}
where $\bd{n}$ is the $n$-th expansion coefficient. Note that the expansion
above is in the moving reference frame and scaled by local temperature. It is
exactly the same expansion used by Grad~\cite{Grad1949a}.  Denoting the
expansion coefficients of $\f{0}$ and $\f{1}$ by $\bd{n}_0$ and $\bd{n}_1$
respectively, we have $\bd{0}_0 = \rho$, $\bd{n}_0 = 0$ for all $n>0$, and:
\begin{equation}
	\f{1} = \omega(\bv)\sum_{n=2}^\infty\frac 1{n!}\bd{n}_1:\bH{n}(\bv).
\end{equation}
The collision operator is defined \textit{via} its Hermite expansion
coefficients as:
\begin{equation}
	\label{eq:mrt}
	\bd{n}_\Omega = -\frac 1{\tau_n}\bd{n}_1,
	\quad\mbox{for}\quad n = 2, \cdots, \infty.
\end{equation}
It can be verified to satisfy Eq.~(\ref{eq:def}).

\subsection{Tensor decomposition}

As $\bd{n}_\Omega$ and $\bd{n}_1$ are tensors consisting of multiple components,
a question arises as if the relaxation can be made at a finer scale. We note
that the tensor space of a given rank can be decomposed into subspaces that
furnish the \textit{irreducible representations} of the rotation group SO(3)
which are the minimal units closed under spatial rotation~\cite{Zee2016}.
Particularly, a rank-$n$ fully symmetric tensor can be decomposed into
irreducible components by repeatedly subtract from it its rank-$(n-2)$
traces~\cite{Jerphagnon1978}.  Assuming Einstein summation convention, the
explicit decomposition of the 2nd, 3rd and 4th rank symmetric tensors in
$d$-dimensions are~\cite{Spencer1970}:
\begin{subequations}
	\label{eq:decomp}
\begin{eqnarray}
	a_{ij}   &=& a'_{ij} + \frac 1da_{pp}\delta_{ij}, \\
	a_{ijk}  &=& a'_{ijk} + \frac 1{d+2}
	\left(a_{ppi}\delta_{jk}+a_{ppj}\delta_{ik}+a_{ppk}\delta_{ij}\right),\\
	a_{ijkl} &=& a'_{ijkl} + \frac 1{d+4}
	\left(a''_{ppij}\delta_{kl} + a''_{ppik}\delta_{jl}
	 	+ a''_{ppil}\delta_{jk} + a''_{ppjk}\delta_{il}
	  	+ a''_{ppjl}\delta_{ik} + a''_{ppkl}\delta_{ij}\right)\nonumber\\
			&+& \frac 1{d(d+2)}a_{ppqq}\left(
				\delta_{ij}\delta_{kl}
			  + \delta_{ik}\delta_{jl}
			  + \delta_{il}\delta_{jl}\right),
\end{eqnarray}
\end{subequations}
where $a'_{ij}$, $a'_{ijk}$ and $a'_{ijkl}$ are all traceless, meaning that
the contractions with respect to any pair of indexes vanishes, and
\begin{equation}
	a''_{ppij}\equiv a_{ppij} - \frac 1d a_{ppqq}\delta_{ij}
\end{equation}
is a traceless second rank tensor. Denote the $k$-th irreducible components of
$\bd{n}_1$ and $\bd{n}_\Omega$ by $\bd{n}_{1,k}$ and $\bd{n}_{\Omega,k}$
respectively in the orders of the terms on the right-hand-side of
Eqs.~(\ref{eq:decomp}). The collision operator defined by Eq.~(\ref{eq:mrt}) can
be further refined as:
\begin{equation}
	\label{eq:domega}
	\bd{n}_{\Omega,k} = -\frac{1}{\tau_{nk}}\bd{n}_{1,k},
\end{equation}
where $k\in \{1, 2\}$ for $n=2, 3$, and $k\in\{1, 2, 3\}$ for $n=4$.  Thus, up
to the 4th order, we can have up to seven independent relaxation times:
$\tau_{21}$ and $\tau_{22}$ for $\bd{2}_\Omega$, $\tau_{31}$ and $\tau_{32}$ for
$\bd{3}_\Omega$, and $\tau_{41}$, $\tau_{42}$, $\tau_{43}$ for $\bd{4}_\Omega$.
Among them, $\tau_{21}$, $\tau_{22}$ and $\tau_{32}$ dictate respectively the
shear viscosity, bulk viscosity in a gas with internal degrees of freedom, and
thermal diffusivity. The other four do not explicitly manifest in the
Navier-Stokes-Fourier equations.

\subsection{Binomial transform}

Direct determination of $\bd{n}_1$ and $\bd{n}_\Omega$ from $f(\bv_i)$ by
Eqs.~(\ref{eq:iso}) requires an interpolation scheme of some kind as $\bv_i$
depends on $\bu$ and $\theta$ and varies with space and
time~\cite{Sun2000,Dorschner2018}.  To avoid the associated numerical
dissipation, the \textit{binomial transform}, Eq.~(\ref{eq:da}), between the
Hermite expansion coefficients in the absolute and relative frame can be used to
obtain $\bd{n}_1$ from $\ba{n}_1$ as:
\begin{subequations}
\label{eq:a2d}
\begin{eqnarray}
	\theta^{\frac{D+2}2}\bd{2}_1 &=& \ba{2}_1,\\
	\theta^{\frac{D+3}2}\bd{3}_1 &=& \ba{3}_1 - 3\bu\ba{2}_1,\\
	\theta^{\frac{D+4}2}\bd{4}_1 &=& \ba{4}_1 - 4\bu\ba{3}_1
	+ 6\left[\bu\bu + (1-\theta)\bdelta\right]\ba{2}_1,
\end{eqnarray}
\end{subequations}
and from $\bd{n}_\Omega$ to $\ba{n}_\Omega$ as:
\begin{subequations}
	\label{eq:d2a}
\begin{eqnarray}
	\ba{2}_\Omega &=& \theta^{\frac{D+2}2}\bd{2}_\Omega,\\
	\ba{3}_\Omega &=& \theta^{\frac{D+3}2}\bd{3}_\Omega + 3\bu\ba{2}_\Omega,\\
	\ba{4}_\Omega &=& \theta^{\frac{D+4}2}\bd{4}_\Omega + 4\bu\ba{3}_\Omega
		 - 6\left[\bu\bu + (1-\theta)\bdelta\right]\ba{2}_\Omega.
\end{eqnarray}
\end{subequations}
Here the fact that $\ba{n}_1 = \ba{n}_\Omega = 0$ for $n = 0, 1$ is used. On
substituting Eqs.~(\ref{eq:a2d}) into Eqs.~(\ref{eq:d2a}), we have:
\begin{subequations}
	\begin{eqnarray}
		\ba{2}_\Omega &=& -\omega_2\ba{2}_1,\\
		\ba{3}_\Omega &=& -\omega_3\ba{3}_1
			+ 3\left(\omega_3-\omega_2\right)\bu\ba{2}_1,\\
		\label{eq:21c}
		\ba{4}_\Omega &=& -\omega_4\ba{4}_1
			+ 4\left(\omega_4-\omega_3\right)\bu\ba{3}_1\nonumber\\
			&&- 6\left[(\omega_4+\omega_2-2\omega_3)\bu\bu
				+ (\omega_4-\omega_2)(1-\theta)\bdelta\right]\ba{2}_1.
	\end{eqnarray}
\end{subequations}
The only difference with the previous model without temperature
scaling~\cite{Shan2019,Li2019} is the additional factor of $1-\theta$ in
Eq.~(\ref{eq:21c}).

\subsection{Space-time discretization}

To numerically solve Eq.~(\ref{eq:lbgk}), the time and spatial derivatives on
the left-hand-side must be descritized.  Integrating using the first-order
forward-Euler scheme from $t=0$ to $1$, we have:
\begin{equation}
	\label{eq:lbgkd}
	f_i(\bx + \bxi_i, t+1) - f_i(\bx, t) = \Omega(f_i).
\end{equation}
It is well-known~\cite{Chen1998a} that the effect of the implicit second-order
error can be absorbed into the dissipation term, effectively making the scheme
second-order accurate with a viscosity proportional to $\tau+1/2$ instead of
$\tau$.  Alternatively, a second-order implicit scheme can be obtained by
integrating Eq.~(\ref{eq:lbgk}) using the trapezoidal rule~\cite{He1998}:
\begin{equation}
	\label{eq:23}
	f_i(\bx + \bxi_i, t+1) - f_i(\bx, t) = \frac 12
	\left[\Omega(f_i(\bx + \bxi_i, t+1)) + \Omega(f_i(\bx, t))\right].
\end{equation}
With the BGK operator, by introducing an auxiliary variable, it was shown that
the dynamics of the distribution at the second order is equivalent to that at
the first order with $\tau$ replaced by $\tau+1/2$~\cite{He1998}. With the more
complicated collision operator defined above, such a change-of-variable is not
directly possible. Nevertheless the dynamics of the moments can be analyzed in a
similar fashion. Define the moment operator $M_{n,k}(f_i)$ as the $k$-th
component of the $n$th order moment of $f_i$. Noting that:
\begin{equation}
	M_{n,k}\left(\Omega(f_i)\right) = -\frac{1}{\tau_{nk}}
	\left[M_{n,k}(f_i)-M_{n,k}(\f{0}_i)\right],
\end{equation}
by taking moments of Eq.~(\ref{eq:23}) and noting $f_i(\bx + \bxi_i, t+1)$ is
the out-going distribution after collision, we write:
\begin{equation}
	M_{n,k}^{out} - M_{n,k}^{in} = -\frac{1}{\tau_{nk}}
	\left[\frac 12\left(M_{n,k}^{out} + M_{n,k}^{in}\right)-M_{n,k}^{eq}\right].
\end{equation}
where the superscripts $out$, $in$ and $eq$ denote respectively the out-going,
in-coming and equilibrium moments.  After some straightforward manipulations, we
have:
\begin{equation}
	M_{n,k}^{out} - M_{n,k}^{in} = -\frac{1}{\tau_{nk}+1/2}
	\left[M_{n,k}^{in}-M_{n,k}^{eq}\right].
\end{equation}
Evidently, at the second order of space-time discretization, the behavior of the
moments, consequently the hydrodynamics, is the same as that at the first order
but with $\tau_{nk}$ replaced by $\tau_{nk}+1/2$. This correspondence is
essentially a consequence of the orthogonal relaxation of the moments.

The actual implementation of the collision operator goes as the following. From
the post-streaming $f_i$, we compute $\f{0}_i$ and then $\f{1}_i$. Using
Eq.~(\ref{eq:iso}) we obtain $\ba{n}_1$ and then $\bd{n}_1$ using
Eqs.~(\ref{eq:a2d}).  Decomposing $\bd{n}_1$ according to Eqs.~(\ref{eq:decomp})
if necessary, applying the relaxation factor  according to
Eqs.~(\ref{eq:domega}) and re-assembling, we obtain $\bd{n}_\Omega$ and then
$\ba{n}_\Omega$ using Eqs.~(\ref{eq:d2a}).  Using the second part of
Eq.~(\ref{eq:iso}) to calculate the bracket on the right-hand-side of
Eq.~(\ref{eq:lbgkd}) all together which effectively trims the part of the
distribution function lying outside of the functional space spanned by the
Hermite polynomials.

\section{Numerical verification}

\label{sec:numerical}

To verify the model, the transport coefficients are measured form the dynamics
of the linear hydrodynamic modes and compared with their theoretical values. The
case setup has been extensively discussed
previously~\cite{Shan2007,Li2011a,Shan2019}.  Here we briefly summarize the
analytical results.  Consider the monochromatic plane wave perturbation:
\begin{equation}
	\left(
	\begin{array}{c}
		\rho \\
		\bu \\
		\theta \\
	\end{array}
	\right)
	=
	\left(
	\begin{array}{c}
		\rho_0 \\
		\bu_0 \\
		\theta_0 \\
	\end{array}
	\right)
	+
	\left(
	\begin{array}{c}
		\bar{\rho} \\
		\bar{\bu} \\
		\bar{\theta} \\
	\end{array}
	\right)
	e^{\omega t + i\bk \cdot (\bm{x}-\bu_0 t) }
\end{equation}
where the subscript $_0$ denotes the base flow and $\bar{\rho}$, $\bar{\bu}$ and
$\bar{\theta}$ perturbation amplitudes. $\bm{k}$ and $\omega$ respectively are
the \textit{wave vector} and \textit{angular frequency} of the plane wave.
Decomposing the velocity into components parallel and perpendicular to the wave
vector and substituting into the Navier-Stokes-Fourier equation, we obtain an
eigen-system in the linear space of $(\bar{\rho}, \bar{u}_\parallel,
\bar{\theta}, \bar{u}_\perp)^T$ where $\bar{u}_\parallel$ and $\bar{u}_\perp$
are the amplitudes of velocity perturbation in the parallel and perpendicular
directions.  The dispersion relations of the four linear modes are:
\begin{subequations}
	\label{eq:dr}
	\begin{eqnarray}
		\omega_v &=& -\nu k^2,\\
		\omega_t &=& -\kappa k^2 \left[1 + \frac{(\gamma-1)\lambda}{\pe^2}\right]
		+ \mathcal{O}\left(\frac{1}{\pe^4}\right),\\
		\label{eq:dr-3}
		\omega_\pm &=& -\alpha k^2
		\left[1-\frac{(\gamma-1)\lambda}{(\gamma-\lambda)\pe^2}\right]
		\pm ic_sk\left[1-\frac{(\gamma +\lambda)^2 - 4\lambda}{8\pe^2}\right]
		+ \mathcal{O}\left(\frac{1}{\pe^4}\right),
	\end{eqnarray}
\end{subequations}
where $\omega_v$, $\omega_t$, and $\omega_\pm$ are the angular frequencies of
the viscous, thermal and two acoustic modes, $c_s\equiv\sqrt{\gamma\theta_0}$ is
a {\em characteristic} speed of sound, $k\equiv |\bm{k}|$ the \textit{wave
	number}, $\re\equiv c_s/\nu k$, $\pe\equiv c_s/\kappa k$ and
$\pr\equiv\nu/\kappa$ the \textit{acoustic} Reynolds, P\'eclet and Prandtl
numbers, and
\begin{equation}
	\lambda\equiv 1 - \left(2-\frac 2D +\frac{\nu_b}{\nu}\right)\pr
\end{equation}
a constant defined for brevity. $\alpha$ is the leading order sound attenuation
rate which is a weighted sum of the shear viscosity, bulk viscosity and thermal
diffusivity:
\begin{equation}
	\alpha = \frac{\gamma-1}2\kappa + \frac{D-1}D\nu + \frac 12\nu_b.
\end{equation}
All the transport coefficients above are given by the relaxation times and the
internal degree of freedom, $S$, by:
\begin{equation}
	\label{eq:trans}
	\nu = \theta_0\left(\tau_{21}-\frac 12\right),\quad
	\nu_b = \frac{2S\theta_0}{D(D+S)}\left(\tau_{22}-\frac 12\right),\quad
	\kappa = \theta_0\left(\tau_{32}-\frac 12\right),\quad
	\gamma = 1 + \frac{2}{D+S}.
\end{equation}

Some remarks are called for at this point. First, while the viscous mode is
independent of the other three and its dispersion relation exact, the dispersion
relations of the thermal and acoustic modes are coupled solutions of a cubic
characteristic equation and only their asymptotic form at large-\pe\ limit are
given above. Second, the bulk viscosity has a leading-order effect on the sound
attenuation but its effects everywhere else is second order. In the following
tests we numerically measure $\omega_v$, $\omega_t$, and $\omega_\pm$ for ranges
of the relaxation times and the internal degree of freedom $S$. Agreement
between the numerical and theoretical values, measured by the relative error:
\begin{equation}
	\omega^\ast\equiv \left|\frac{\omega_{numerical} - \omega_{theoretical}}
		{\omega_{theoretical}}\right|,
\end{equation}
sufficiently verifies the accuracy of the model.

\begin{figure}
	\begin{subfigure}{.5\textwidth}
		\centering
		\includegraphics[width=\linewidth]{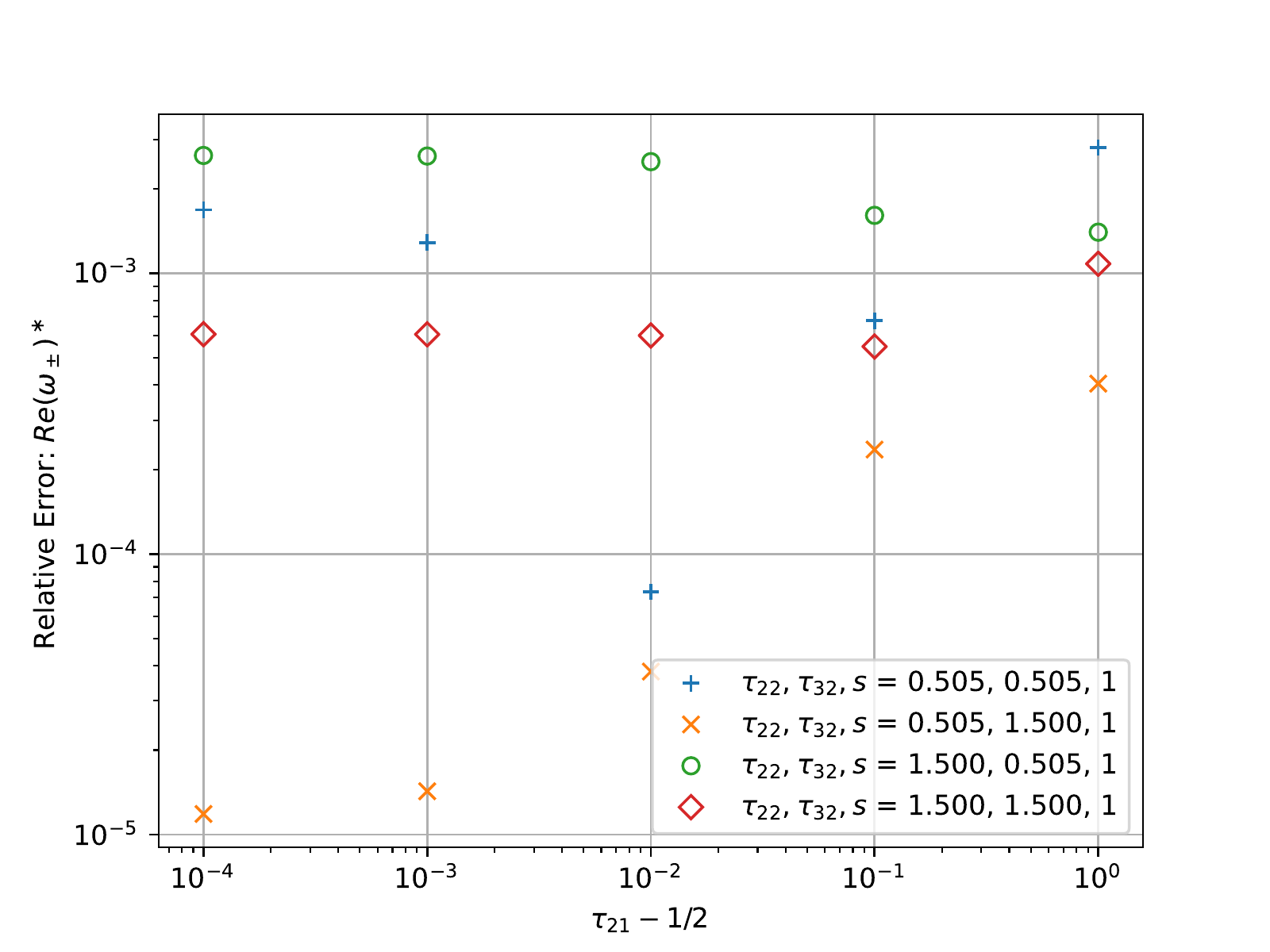}
	\end{subfigure}%
	\begin{subfigure}{.5\textwidth}
		\centering
		\includegraphics[width=\linewidth]{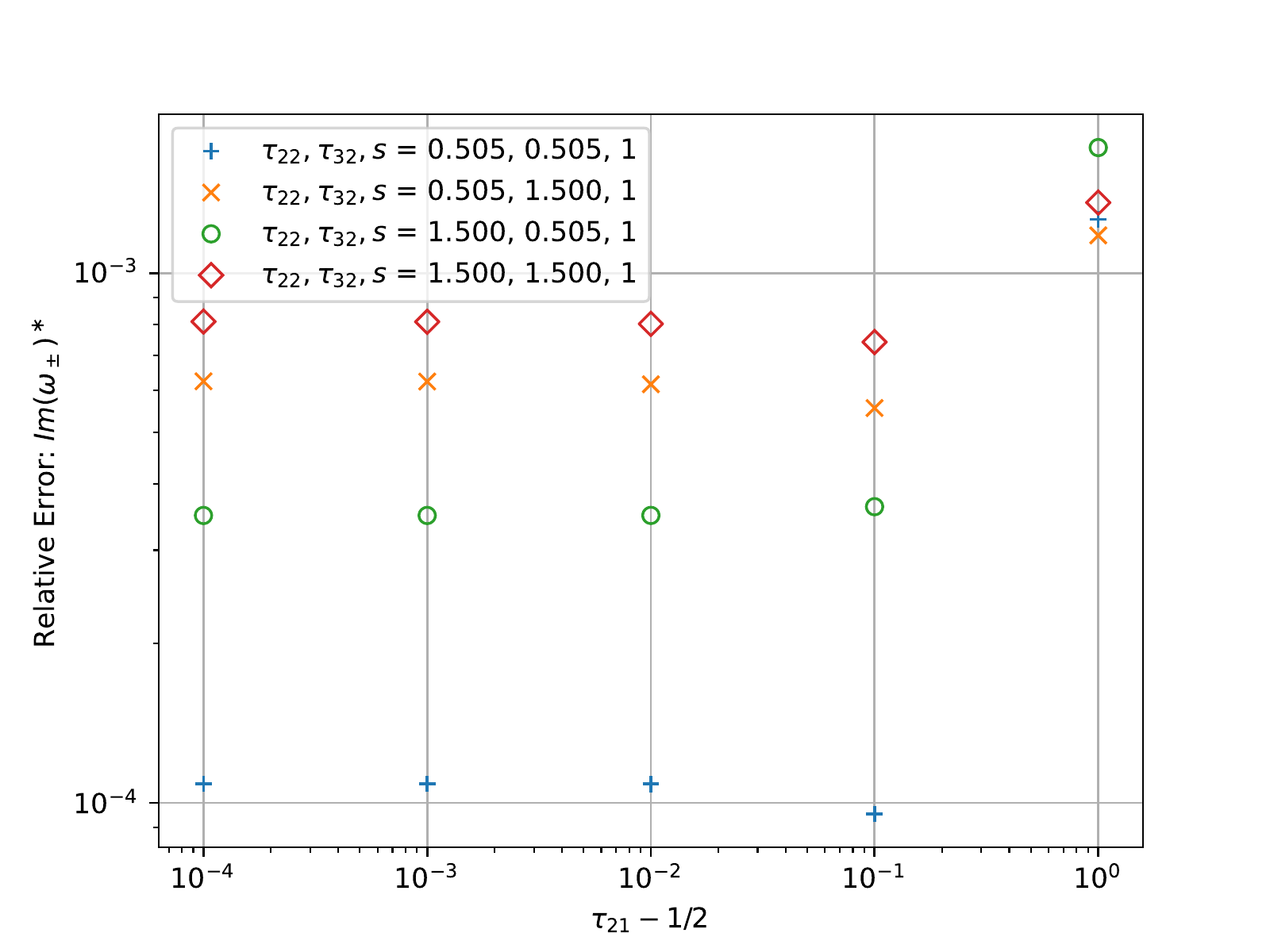}
	\end{subfigure}
	\begin{subfigure}{.5\textwidth}
		\centering
		\includegraphics[width=\linewidth]{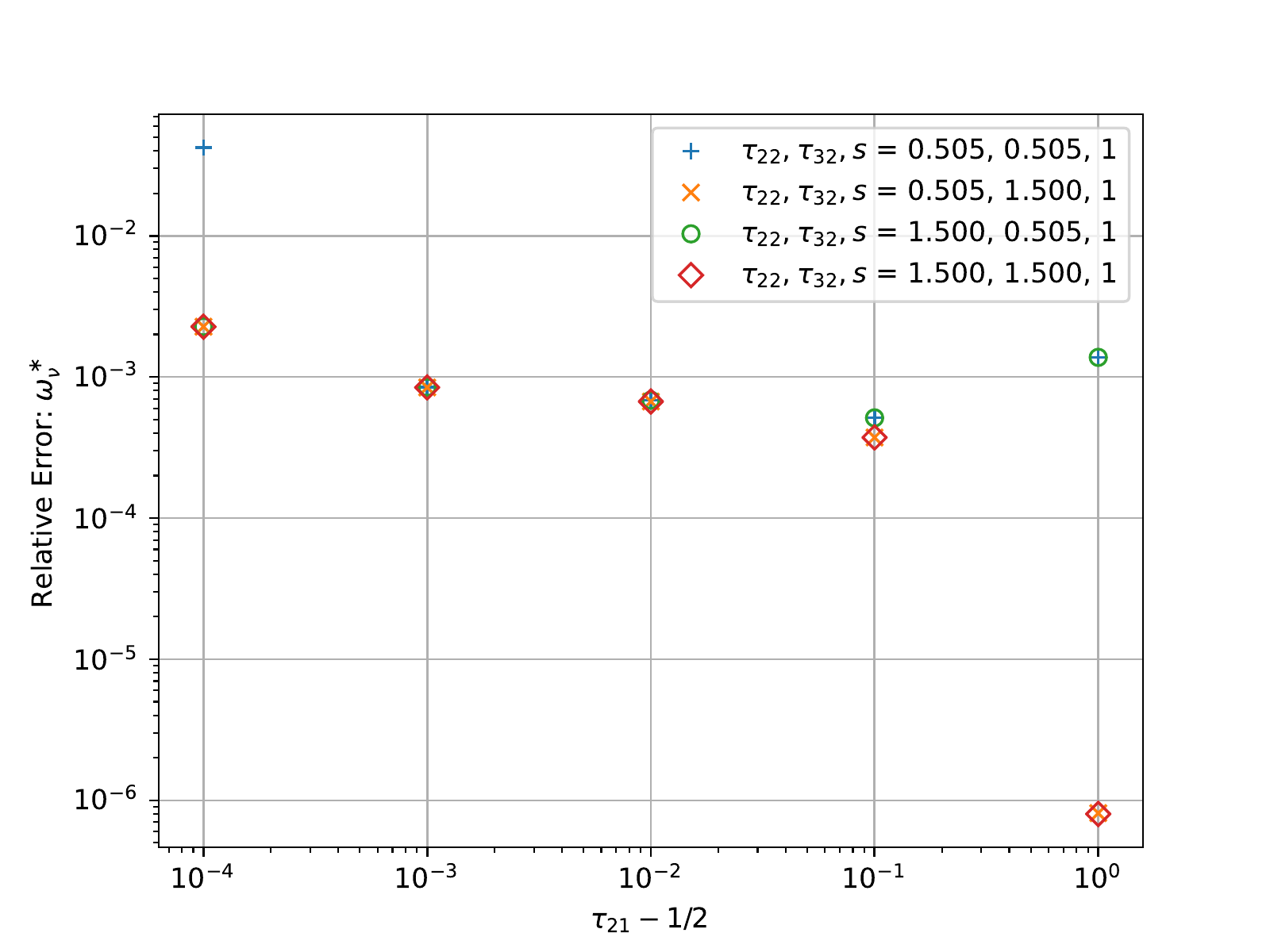}
	\end{subfigure}%
	\begin{subfigure}{.5\textwidth}
		\centering
		\includegraphics[width=\linewidth]{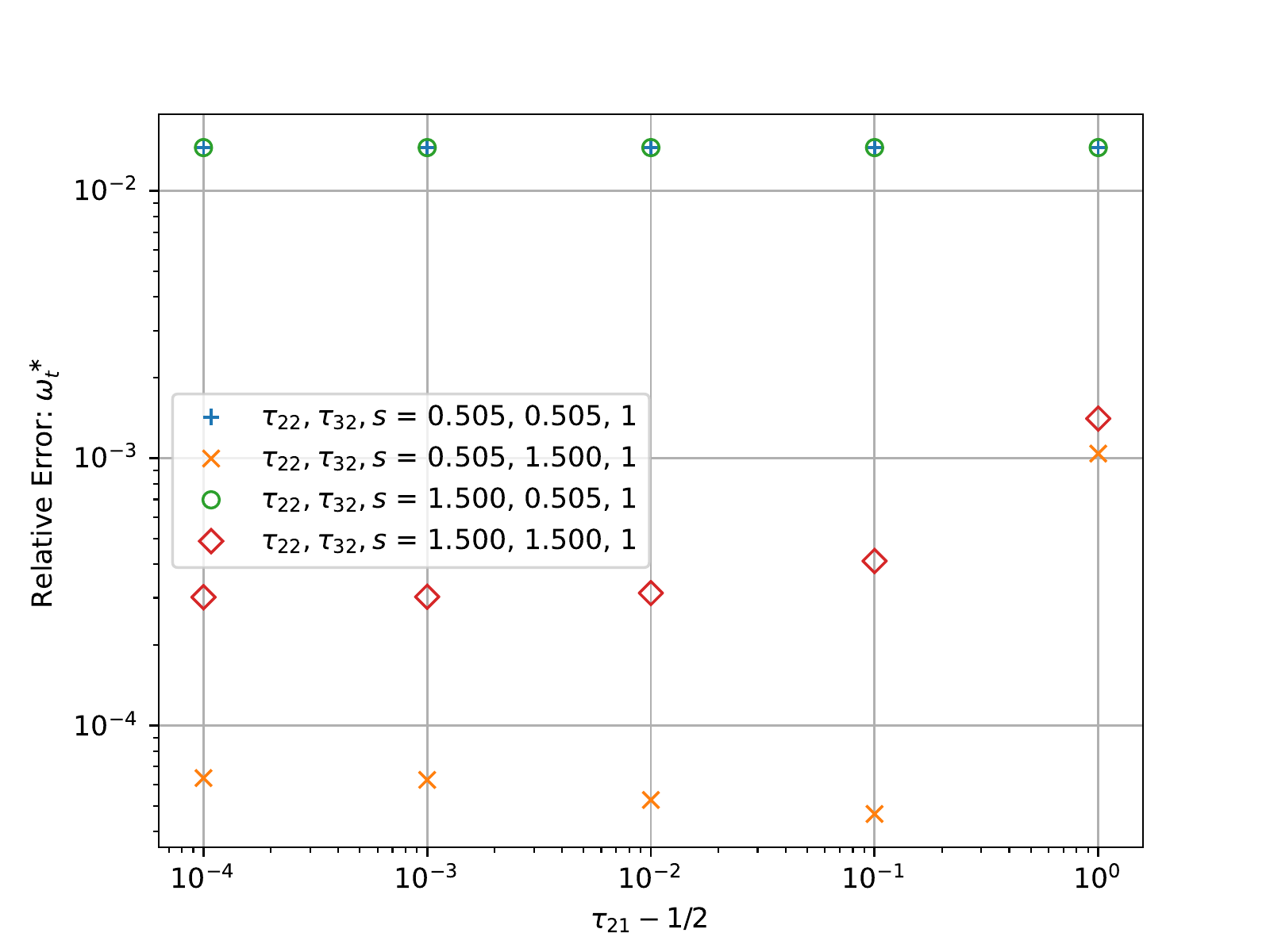}
	\end{subfigure}
	\caption{The relative errors in the four eigen-frequencies of the hydrodynamic
	modes for a wide range of the three relaxation times $\tau_{21}$, $\tau_{22}$
	and $\tau_{32}$.}

	\label{fig:fig}
\end{figure}

Shown in Fig.~\ref{fig:fig} are the relative errors in $\omega_v$, $\omega_t$,
and $\omega_\pm$ for ranges of the relaxation times and the internal degree of
freedom $S$.  The angular frequencies are extracted by fitting the time history
of the amplitudes of the four modes.  The numerical simulations were conducted
in the same fashion reported previously~\cite{Li2011a,Shan2019} on a lattice of
$100\times 100$.  It can be seen that the relative error in the angular
frequencies are generally less than 1\%, verifying the correctness of
Eqs.~(\ref{eq:trans}).

\section{Conclusions and discussions}

\label{sec:conclusion}

In the present work we present in detail a generic multiple relaxation collision
model based on Hermite expansion and tensor decomposition.  The collision
operator is first expanded in terms of the Hermite tensorial polynomials in the
reference frame moving with the fluid.  Each of the polynomials is further
decomposed into traceless components corresponding to the irreducible
representations of the rotation group and relaxed with independent relaxation
rate which is the most generic form of relaxation without violating rotation
symmetry. The number of the independent relaxation rates are given by the number
of traceless tensors of each tensor rank and is 2, 2, and 3 for the 2nd, 3rd,
and 4th rank tensors respectively.  To avoid possible numerical dissipation due
to interpolation, the expansion coefficients are transformed into laboratory
frame for evaluation with fixed-velocity quadrature rules.  In its crudest form
the model reverts to the BGK model.  The first few relaxation rates are shown to
correspond to the shear and bulk viscosity and the thermal diffusivity of the
translational and other forms of motion.  Numerical measurements of the
eigen-frequencies of the linear hydrodynamic modes agree well with theoretical
predictions. Although the other relaxation times do not appear in the
hydrodynamic equations and have negligible effect in continuum, their effects
are presumably more explicit in situations such as rarefied gases where the
distribution is away from the Maxwell-Boltzmann equilibrium.  At least in the
linear regime, the spectrum of the relaxation times dictates the behavior of the
collision operator and could potentially serve as a set of definitive variables.

This work was supported by National Science Foundation of China Grants~91741101
and 91752204, Department of Science and Technology of Guangdong Province Grant
2019B21203001, Shenzhen Science and Technology Program Grant
KQTD20180411143441009.

\appendix

\section{Transforms between the central and raw moments}

The $n$-th order Hermite polynomial in $D$-dimensions is defined by the
Rodrigues' formula:
\begin{equation}
	\label{eq:rodrigues}
	\bH{n}(\bxi) = \frac{(-1)^n}{\omega(\bxi)}\nabla^n\omega(\bxi),
\end{equation}
where $\omega(\bxi)$ is the weight function:
\begin{equation}
	\label{eq:wf}
	\omega(\bxi) = \frac 1{(2\pi)^{D/2}}e^{-\xi^2/2},
\end{equation}
with $\xi^2 = \bxi\cdot\bxi$.  $\bH{n}(\bxi)$ is a fully symmetric rank-$n$
tensor and an $n$-th degree polynomial in $\bxi$. The explicit form of
$\bH{n}(\bxi)$ is:
\begin{equation}
	\label{eq:hermite-explicit}
	\bH{n}(\bxi) = \sum_{k=0}^{\lfloor n/2\rfloor}(-1)^kD_n^k\bxi^{n-2k}\bdelta^k,
\end{equation}
where the floor function, $\lfloor x\rfloor$, stands for the largest integer not
exceeding $x$, and:
\begin{equation}
	D_n^k \equiv\frac{n!}{(n-2k)!2^kk!} = C_n^{2k}(2k-1)!!,
\end{equation}
where $C_n^k$ is the \textit{binomial coefficient}. The leading values of
$D_n^k$ are listed in Tab.~\ref{tb:dnk} for convenience.
\begin{table}[htb]
	\centering
	\begin{tabular}{rr|rrrrrrrr}
		\hline\hline
		                     &   &      \multicolumn{8}{c}{$n$}      \\
		                     &   & 0 & 1 & 2 & 3 & 4 &  5 &  6 &   7 \\ \hline
		\multirow{4}{*}{$k$} & 0 & 1 & 1 & 1 & 1 & 1 &  1 &  1 &   1 \\
		                     & 1 &   &   & 1 & 3 & 6 & 10 & 15 &  21 \\
		                     & 2 &   &   &   &   & 3 & 15 & 45 & 105 \\
		                     & 3 &   &   &   &   &   &    & 15 & 105 \\ \hline\hline
	\end{tabular}
	\caption{Leading values of the coefficients $D_n^k$.}
	\label{tb:dnk}
\end{table}

Differentiating Eqs.~(\ref{eq:rodrigues}) and (\ref{eq:hermite-explicit}) and
noticing $\nabla\omega^{-1} = \bxi/\omega$, we have respectively:
\begin{equation}
	\nabla\bH{n}(\bxi) = (-1)^n\left(\frac{\nabla^{n+1}\omega}\omega +
	\frac{\bxi\nabla^n\omega}\omega\right) = -\bH{n+1}(\bxi) + \bxi\bH{n}(\bxi),
\end{equation}
and:
\begin{equation}
	\nabla\bH{n}(\bxi) = \sum_{k=0}^{\lfloor n/2\rfloor}
	(n-2k)D_n^k\bxi^{n-1-2k}\bdelta\bdelta^k = n\bdelta\bH{n-1}(\bxi).
\end{equation}
Combining the two equations above, we have the recurrence relation:
\begin{equation}
	\label{eq:recur}
	\bxi\bH{n}(\bxi) = \bH{n+1}(\bxi) + n\bdelta\bH{n-1}(\bxi).
\end{equation}

We first establish by induction that the monomials can be expressed by Hermite
polynomials as the following:
\begin{equation}
	\label{eq:monomials}
	\bxi^n = \sum_{k=0}^{\lfloor n/2\rfloor}D_n^k\bH{n-2k}(\bxi)\bdelta^k.
\end{equation}
It is trivially true for $n=0$. Multiple the above by $\bxi$ and using
Eq.~(\ref{eq:recur}), we have:
\begin{eqnarray}
	\lefteqn{\bxi^{n+1} = \sum_{k=0}^{\lfloor n/2\rfloor}D_n^k
		\left[\bH{n+1-2k}(\bxi)+(n-2k)\bdelta\bH{n-1-2k}(\bxi)\right]\bdelta^k}\nonumber\\
	 &=& \sum_{k=0}^{\lfloor n/2\rfloor}D_n^k\bH{n+1-2k}(\bxi)\bdelta^k
+\sum_{k=0}^{\lfloor n/2\rfloor}(n-2k)D_n^k\bH{n-1-2k}(\bxi)\bdelta^{k+1}.
\end{eqnarray}
Applying the change of variable $k+1\rightarrow k$ to the second term, we have:
\begin{eqnarray}
	\bxi^{n+1} &=& \sum_{s=0}^{\lfloor (n+1)/2\rfloor}
	\left[1+\frac{2k}{(n+1-2k)}\right]D_n^k
	\bH{n+1-2k}(\bxi)\bdelta^k\nonumber\\
	&=& \sum_{k=0}^{\lfloor (n+1)/2\rfloor}D_{n+1}^k
		\bH{n+1-2k}(\bxi)\bdelta^k,
\end{eqnarray}
completing the induction.

We now establish the following relation by induction:
\begin{equation}
	\label{eq:hc}
	\bH{n}(\bxi+\bu) = \sum_{k=0}^nC_n^k\bH{k}(\bxi)\bu^{n-k},
\end{equation}
Again it is trivially true for $n = 0$. Using Eq.~(\ref{eq:recur}) repeatedly,
we have
\begin{eqnarray}
	\lefteqn{\bH{n+1}(\bxi+\bu) = (\bxi+\bu)\bH{n}(\bxi+\bu) -
		n\bdelta\bH{n-1}(\bxi+\bu)}\nonumber\\
	&=&(\bxi+\bu)\sum_{k=0}^nC_n^k\bH{k}(\bxi)\bu^{n-k}
	- n\bdelta\sum_{k=0}^{n-1}C_{n-1}^k\bH{k}(\bxi)\bu^{n-1-k}\nonumber\\
	&=&\sum_{i=0}^nC_n^k \bH{k+1}(\bxi)\bu^{n-k}
	+\bdelta\sum_{k=0}^nkC_n^k\bH{k-1}(\bxi)\bu^{n-k}\nonumber\\
	&&+ \sum_{k=0}^nC_n^k\bH{k}(\bxi)\bu^{n+1-k}
	  -n\bdelta\sum_{k=0}^{n-1}C_{n-1}^k\bH{k}(\bxi)\bu^{n-1-k}.
\end{eqnarray}
Noticing that $nC_{n-1}^k = (k+1)C_n^{k+1}$, after a change of variable
$k+1\rightarrow k$, the second term cancels the last term, and the first term
becomes:
\begin{equation}
	\sum_{k=1}^{n+1}C_n^{k-1}\bH{k}(\bxi)\bu^{n+1-k}.
\end{equation}
On combining with the third term and noting that if we define $C_n^{-1} \equiv
C_n^{n+1}\equiv 0$,
\begin{equation}
	C_n^{k-1} + C_n^k = C_{n+1}^k,\quad\mbox{for}\quad k = 0, \cdots, n,
\end{equation}
we have:
\begin{equation}
	\bH{n+1}(\bxi+\bu) = \sum_{k=0}^{n+1}C_{n+1}^k\bH{k}(\bxi)\bu^{n+1-k},
\end{equation}
completing the induction.

Now consider the Hermite polynomials under coordinate scaling by a constant
factor, $\alpha$.  Using Eq.~(\ref{eq:hermite-explicit}), we have:
\begin{equation}
	\bH{n}(\alpha\bxi) = \sum_{k=0}^{\lfloor n/2\rfloor}
	(-1)^k\alpha^{n-2k}D_n^k\bxi^{n-2k}\bdelta^k.
\end{equation}
On substituting Eq.~(\ref{eq:monomials}) into the equation above, we have:
\begin{eqnarray}
	\bH{n}(\alpha\bxi) &=& \sum_{k=0}^{\lfloor n/2\rfloor}
	(-1)^k\alpha^{n-2k}D_n^k\bdelta^k
	\left[\sum_{t=0}^{\lfloor n/2-k\rfloor}
	D_{n-2k}^t\bH{n-2k-2t}(\bxi)\bdelta^t\right]\nonumber\\
	&=&\sum_{k=0}^{\lfloor n/2\rfloor}\sum_{t=0}^{\lfloor n/2-k\rfloor}
	\frac{(-1)^kn!\alpha^{n-2k}}{(n-2k-2t)!2^{k+t}k!t!}\bH{n-2k-2t}(\bxi)\bdelta^{k+t}.
\end{eqnarray}
Realizing that the double summation is over all combinations of $k$ and $t$ such
that $k+t \leq\lfloor n/2\rfloor$, we define $m=k+t$ and re-arrange the
summation to write:
\begin{eqnarray}
	\bH{n}(\alpha\bxi) &=& \alpha^n\sum_{m=0}^{\lfloor n/2\rfloor}\frac{n!}{(n-2m)!2^mm!}
	\left[\sum_{k=0}^m\frac{m!(-\alpha^{-2})^k}{k!(m-k)!}\right]
	\bH{n-2m}(\bxi)\bdelta^m\nonumber\\
	\label{eq:hermite-scale}
	&=& \alpha^n\sum_{m=0}^{\lfloor n/2\rfloor}(1-\alpha^{-2})^mD_n^m
	\bH{n-2m}(\bxi)\bdelta^m.
\end{eqnarray}

Combining Eqs.~(\ref{eq:hc}) and (\ref{eq:hermite-scale}), we can express the
Hermite polynomials in the moving reference frame scaled by temperature using
the ones in the absolute frame. We write:
\begin{eqnarray}
	\lefteqn{\theta^{\frac n2}\bH{n}\left(\frac{\bxi-\bu}{\sqrt{\theta}}\right) =
		\sum_{k=0}^{\lfloor n/2\rfloor}D_n^k(1-\theta)^k
		\bH{n-2k}(\bxi-\bu)\bdelta^k}\nonumber\\
	&=& \sum_{k=0}^{\lfloor n/2\rfloor}D_n^k(1-\theta)^k
	\left[\sum_{m=0}^{n-2k}C_{n-2k}^m\bH{m}(\bxi)(-\bu)^{n-2k-m}\right]\bdelta^k\nonumber\\
	&=& \sum_{k=0}^{\lfloor n/2\rfloor}\sum_{m=0}^{n-2k}
	\frac{(-1)^{n-m-2k}n!(1-\theta)^k}{(n-m-2k)!k!m!2^k}
	\bH{m}(\bxi)\bu^{n-m-2k}\bdelta^k.\nonumber
\end{eqnarray}
Again, realizing that the summation is over all $k$ and $m$ such that
$m+2k\leq n$, the order of the two summations can be swapped:
\begin{equation}
	\sum_{k=0}^{\lfloor n/2\rfloor}\sum_{m=0}^{n-2k} =
	\sum_{m=0}^n\sum_{k=0}^{\lfloor (n-m)/2\rfloor}.
\end{equation}
Also notice that $(-1)^{n-m-2k} = (-1)^{n-m}$ and
\begin{equation}
	\frac{n!}{(n-m-2k)!k!m!2^k} = \frac{n!(n-m)!}{(n-m)!m!(n-m-2k)!k!2^k}
	 = C_n^mD_{n-m}^k,
\end{equation}
we have:
\begin{equation}
	\theta^{\frac n2}\bH{n}\left(\frac{\bxi-\bu}{\sqrt{\theta}}\right) =
	\sum_{m=0}^n(-1)^{n-m}C_n^m\bH{m}(\bxi)
	\sum_{k=0}^{\lfloor (n-m)/2\rfloor}D_{n-m}^k(1-\theta)^k
	\bu^{n-m-2k}\bdelta^k,
\end{equation}
where the inner summation is a function of $\bu$ and $\theta$.  Defining:
\begin{equation}
	\bm{A}_m(\bu, \theta) = \sum_{k=0}^{\lfloor m/2\rfloor}
	D_m^k(1-\theta)^k\bu^{m-2k}\bdelta^k,
	\quad\mbox{with}\quad
	\bm{A}_0 = 1,\quad \bm{A}_1 = \bu,
\end{equation}
we can finally write:
\begin{equation}
	\label{eq:hv}
	\bH{n}\left(\frac{\bxi-\bu}{\sqrt{\theta}}\right) = \theta^{-\frac n2}
	\sum_{k=0}^n(-1)^{n-k}C_n^k\bH{k}(\bxi)\bm{A}_{n-k}(\bu, \theta).
\end{equation}
Note that if $\theta=1$, $\bm{A}_m = \bu^m$ and the above reverts to
Eq.~(\ref{eq:hc}).

Denoting $\bv = (\bxi-\bu)/\sqrt{\theta}$, Hermite expansion can be shifted
and/or scaled as:
\begin{subequations}
	\begin{eqnarray}
		f(\bv) &=& \omega(\bv)\sum_{n=0}^\infty\frac 1{n!}\bd{n}:\bH{n}(\bv),\\
		\label{eq:expansion-v}
		\bd{n} &=& \int f(\bv)\bH{n}(\bv)d\bv,\quad n = 0, \cdots, \infty,
	\end{eqnarray}
\end{subequations}
On substituting Eq.~(\ref{eq:hv}) and $d\bv = \theta^{-D/2}d\bxi$ into
Eq.~(\ref{eq:expansion-v}), we have:
\begin{equation}
	\label{eq:da}
	\bd{n} = \theta^{-\frac{D+n}2}\sum_{k=0}^n(-1)^{n-k}C_n^k\ba{k}\bm{A}_{n-k}.
\end{equation}

Noticing that the Maxwell-Boltzmann distribution function is:
\begin{equation}
	\f{0} \equiv \frac{\rho}{(2\pi\theta)^{D/2}}
	\exp\left[-\frac{(\bxi-\bu)^2}{2\theta}\right]
	= \rho\theta^{-\frac D2}\omega(\bv),
\end{equation}
its Hermite expansion coefficients in the absolute
frame~\cite{Shan1998,Shim2012} is:
\begin{equation}
	\ba{n} = \rho\theta^{-\frac D2}\int\omega(\bv)\bH{n}(\bxi)d\bxi
	=\rho\int\omega(\bv)\bH{n}(\sqrt{\theta}\bv + \bu)d\bv.
\end{equation}
Using Eq.~(\ref{eq:hv}) to expand $\bH{n}(\sqrt{\theta}\bv + \bu)$ in terms of
$\bH{n}(\bv)$, owing to the orthogonal relation, only the zero-th order term
survives the integration. We have:
\begin{equation}
	\ba{n} = \rho\theta^{-n/2} (-1)^n\bm{A}_n
	\left(-\sqrt{\theta}\bu, \frac 1\theta\right)
	= \rho\sum_{k=0}^{\lfloor m/2\rfloor}
	D_m^k(\theta-1)^k\bu^{m-2k}\bdelta^k.
\end{equation}


%

\end{document}